\documentclass[12pt]{amsart}
\usepackage{amssymb}
\usepackage{amsbsy}
\usepackage{amscd}
\usepackage{eepic}
\usepackage[dvips]{graphics}
\DeclareGraphicsExtensions{.eps}
\makeatletter
\def\cal{\mathcal}
\def\Bbb{\mathbb}

\newenvironment{pf*}[1]{\proof[#1]}{\endproof}
\renewcommand{\thesubsection}{\thesection(\@roman\c@subsection)}
\makeatother
\newtheorem{Theorem}[equation]{Theorem}
\newtheorem{Corollary}[equation]{Corollary}

\theoremstyle{definition}
\newtheorem{Definition}[equation]{Definition}

\theoremstyle{remark}
\newtheorem{Remark}[equation]{Remark}

\numberwithin{equation}{section}
\numberwithin{figure}{section}

\begin{document}
\title{Special McKay correspondence}
\author{Yukari Ito}
\address{Department of Mathematics, Tokyo Metropolitan University, Hachioji, 
Tokyo 192-0397, Japan}
\email{yukari@comp.metro-u.ac.jp}
\dedicatory{Dedicated to Professor Riemenschneider for his 60th birthday}
\thanks{The author is partially supported by the Grant-in-aid for 
Scientific Research (No.13740019), the Ministry of Education.}
\subjclass{14C05, 14E15}
\maketitle
\begin{abstract}
There are many generalizations of the McKay correspondence for higher
dimensional Gorenstein quotient singularities and there are some
applications to compute the topological invariants today. But some of
the invariants are completely different from the classical invariants,
in particular for non-Gorenstein cases.  
In this paper, we would like to discuss the McKay correspondence for
2-dimensional quotient singularities via \lq\lq special"
representations which gives the classical topological invariants and
give a new characterization of the special representations for cyclic
quotient singularities  in terms of  combinatorics.
\end{abstract}
\tableofcontents
\section{McKay correspondence}

The McKay correspondence is originally a correspondence between the
topology of the minimal resolution of a 2-dimensional rational double
point, which is a quotient singularity by a finite subgroup $G$ of
$SL(2,\Bbb C )$, and the representation theory (irreducible
representations or conjugacy classes) of the group $G$. 
We can see the correspondence via Dynkin diagrams, which came from
McKay's observation in 1979 (\cite{McKay}). 

Let $G$ be a finite subgroup of $SL(2,\Bbb C )$, then the quotient
space  $X:=\Bbb C^2/G$ has a rational double point at the origin.  
As there exists the minimal resolution $\widetilde X$ of the
singularity, we have the exceptional divisors $E_i$. The dual graph of
the configuration of the exceptional divisors is just the Dynkin
diagram of type $A_n$, $D_n$, $E_6$, $E_7$ or $E_8$. 

On the other hand, we have the set of the irreducible representations
$\rho_i$ of the group $G$ up to isomorphism and let $\rho$ be the
natural  representation in $SL(2,\Bbb C )$. The tensor product of these
representations  
$$ \rho_i \otimes \rho = \bigoplus_{j=0}^r a_{ij} \rho_j,$$
where $\rho_0$ is the trivial representation and $r$ is the number of the non-trivial irreducible representations,
gives a set of integers $a_{ij}$ and it determines the Cartan matrix
which defines  the Dynkin diagram. 
\footnote{More precisely, the Cartan matrix is defined as the matrix
  $2E-A$, where  $E$ is the $r\times r$ identity matrix and
  $A=\{a_{ij}\}$ $(i,j \not= 0)$.}  

Then we have a one-to-one numerical correspondence between non-trivial
irreducible representations $\{\rho_i\}$ and irreducible exceptional
curves $\{E_i\}$, that is, the intersection matrix of the exceptional
divisors is the opposite of the Cartan matrix.

This phenomenon was explained geometrically in terms of vector bundles
on the minimal resolution by Gonzalez-Sprinberg and Verdier
(\cite{GV}) \footnote{They gave  the name  {\it McKay correspondence}
  (in French, {\it la correspondance de McKay}) in this paper!} by
case-by-case computations in 1983. In 1985,  
Artin and Verdier (\cite{Artin}) proved this more generally with
reflexive modules and this theory was developed by Esnault and
Kn\"orrer (\cite{EE}, \cite{EK}) for more general quotient surface
singularities. After  Wunram (\cite{Wunram}) constructed a nice
generalized McKay correspondence for any quotient surface
singularities in 1986 in his dissertation,  
Riemenschneider introduced the notion of \lq\lq special
representation etc." and made propaganda for the more generalized
McKay correspondence (cf. \cite{RW}). 
\footnote{Similar generalization for $G\subset GL(2,\Bbb C )$ was
  obtained by Gonzalez-Sprinberg and the related topics were discussed in
   \cite{GS}.} 

In dimension three, we have several \lq\lq McKay correspondences" but
they are just bijections between two sets: 
Let $X$ be  the quotient singularity $\Bbb C^3/G$ where $G$ is a
finite subgroup of $SL(3,\Bbb C)$. Then $X$ has a Gorenstein canonical
singularity of index 1 but not a terminal singularity. It is known
that there exist crepant resolutions  
$\widetilde X$ of this singularity. The crepant resolution is a
minimal resolution and preserves the triviality of the canonical
bundle in this case. 
 
As for the McKay correspondence, the followings are known:
\begin{enumerate}
\item  (Ito-Reid \cite{IR}) There exists a base of cohomology group $H^{2i}(\widetilde X,   \Bbb Q)$, indexed by the conjugacy classes of \lq\lq age"
  $i$ in $G$. 
 
\item  (Ito-Nakajima \cite{YH}) There exists a base of Grothendieck group $K(\widetilde X)$, indexed by the irreducible representations of $G$, when $G$ is a finite abelian subgroup.

\item  (Bridgeland-King-Reid \cite{BKR}) There exists an equivalence between 
the derived category $D(\widetilde X)$ and the equivariant dereived category 
$D^G(\Bbb C^3)$ for any finite subgroup. 

\end{enumerate}

\def\e{\epsilon}

\begin{Remark}
In (1), the {\bf age} of $g\in G$ is defined as follows:
After diagonalization, if $g^r=1$, we obtain 
$g'=\text{diag}(\e^a,\e^b, \e^c)$ where $\e$ is a primitive $r$-th
root of unity. Then $\text{age}(g):=(a+b+c)/r$. For the identity
element {\it  id}, we define age({\it id})$=0$ and all ages are integers if
$G\subset SL(3,\Bbb C)$. 

The correspondence (2) can be included in (3), but note that the
2-dimensional numerical McKay correspondence can be explained very
clearly as a corollary of the result in \cite{YH}.
\end{Remark}

As a generalization of the first McKay correspondence (1), we have a
precise correspondence for each $2i$-th cohomology with conjugacy
classes of age $i$ for any $i=1, \cdots, n-1$ in dimension $n$  which
was given by Batyrev and Kontsevich via \lq\lq motivic integral"
under the assumption of the existence of a crepant resolution, and
this idea was developed to \lq\lq string theoretic cohomology" for all
quotient singularities (cf. \cite{Vitya}).  

 And we can see that the string theoretic Euler number of the
 resolution is the same as the order of the acting group $G$ in case
 $G\subset GL(n,\Bbb C)$,  but it is different from the usual
 topological Euler number of the minimal resolution. 
 Of course, it is very interesting to consider the geometrical meaning
 of these new invariants.

 By the way, in (2) we don't have such a difference among
   representations as {\it age}.   But the  author is interested in
   the relation between the group theory and the classical topological
   invariants.  
   Then we would like to remind the reader of the notion of special
   representations which gives some differences between irreducible
   representations.   
The special representations were defined by Riemenschneider and Wunram
   (\cite{RW}); each of the special irreducible respresentations corresponds 
    to an exceptional divisor of the
   minimal resolution of a 2-dimensional quotient singularity.  

In particular, we would like to discuss special representations and
the minimal resolution for quotient surface singularities from now
on. 
Around 1996, Nakamura and the author showed another way to the McKay
correspondence with the help of the $G$-Hilbert scheme, which is a
2-dimensional $G$-fixed set of the usual Hilbert scheme of $|G|$-points
on $\Bbb C^2$ and isomorphic to the minimal resolution. Kidoh
(\cite{Kidoh})  proved that the $G$-Hilbert scheme for  general cyclic
surface singularities is the minimal resolution. Then Riemenschneider
checked the cyclic case and conjectured that the representations which
are given by the Ito-Nakamura type McKay correspondence via
$G$-Hilbert scheme are just  special representations in 1999
(\cite{Oswald}) and this conjecture was proved by A.~Ishii recently
(\cite{Akira}). 
In this paper, we will give another characterization of  the special
representations by combinatorics for the cyclic quotient case, using
results on the $G$-Hilbert schemes. 

As a colorful introduction to the McKay correspondence, the author
would like to recommend a paper presented at the Bourbaki seminar by
Reid (\cite{Reid})  and also on the Web page
(http://www.maths.warwick.ac.uk/ $\sim$miles/McKay), one can find some
recent papers related to the McKay correspondence. 

This paper is organized as follows:
In this section, we already gave a brief history of the McKay
correspondence and we will discuss the special representations and the
generalized McKay correspondence in the following section. In section
three, we treat $G$-Hilbert schemes as a resolution of singularities,
consider the relation with the toric resolution in the cyclic case,
and show how to find the special representations by combinatorics. In
the final section, we will discuss an example and related topics.

Most of the contents of this paper are based on the author's talk in
the summer school on toric geometry at Fourier Institute in Grenoble,
France in July 2000, and she would like to thank the organizers for
their hospitarity and the participants for the nice atmosphere. She
would like to express her gratitude to Professor Riemenschneider for
giving her a chance to consider the special representations via
$G$-Hilbert schemes and for the various comments and
useful suggestions on her first draft.

\section{Special representations} 

In this section, we will discuss the special representations. Let $G$
 be a finite small subgroup of $GL(2,\Bbb C)$, that is, the action of
 the group  $G$ is free outside the origin,  and  
 $\rho$ be a representation of $G$ on $V$. $G$ acts on 
 $\Bbb C^2 \times V$ and the quotient is a vector bundle on 
 $(\Bbb C^2\setminus\{0\})/G$ 
 which can be extended to a reflexive sheaf $\cal F$ on $X\colon=\Bbb C^2/G$.

For any reflexive sheaf $\cal F$ on a rational surface singularity $X$
and the minimal resolution $\pi\colon\widetilde X \rightarrow X$,
we define a sheaf $\widetilde{\cal F} \colon= \pi^*\cal F /\text{torsion}$.

\def\TF{\widetilde{\cal F}}

\begin{Definition}(\cite{EE})
The sheaf $\TF$ is called a {\it full sheaf} on $\widetilde X$.
\end{Definition}

\begin{Theorem}{\rm (\cite{EE})}
A sheaf $\TF$ on $\widetilde X$ is a full sheaf if the following
conditions are fulfilled: 

1. $\TF$ is locally free,

2. $\TF$ is generated by global sections,   

3. $H^1(\widetilde X, \TF^\vee\otimes \omega_{\widetilde X})=0$, where $\vee$
   means the dual. 
\end{Theorem}

Note that a sheaf $\TF$ is indecomposable if and only if the
corresponding representation $\rho$ is irreducible. Therefore we
obtain an indecomposable full sheaf $\TF_i$ on $\widetilde X$ for each
irreducible representation $\rho_i$, but in general, the number of the
irreducible representations is larger than that of irreducible
exceptional components. Therefore Wunram and Riemenschneider introduced
the notion of speciality for full sheaves: 

\begin{Definition}(\cite{RW})
A full sheaf is called {\it special} if and only if 
$$H^1(\widetilde X, \TF^\vee)=0.$$ 

A reflexive sheaf $\cal F$ on $X$ is {\it special} if $\TF$ is so.

A representation $\rho$ is {\it special} if the associated reflexive
sheaf $\cal F$ on $X$ is special.

\end{Definition}

With these definitions, the following equivalent conditions for the
speciality hold:

\begin{Theorem}\label{Th:RW}{\rm (\cite{RW}, \cite{Wunram})}

1. $\TF$ is special $\Longleftrightarrow$ 
   $\TF\otimes \omega_{\widetilde X}\rightarrow 
   [(\cal F \otimes \omega_{\widetilde X})^{\vee\vee} ]^{\sim}$ 
   is an isomorphism, 

2. $\cal F$ is special $\Longleftrightarrow$ 
   $\cal F \otimes \omega_{\widetilde X}/\text{torsion}$ is reflexive,

3. $\rho$ is a special representation $\Longleftrightarrow$ the map
   ${(\Omega^2_{\Bbb C^2})}^G \otimes (\cal O_{\Bbb C^2} \otimes V)^G 
   \rightarrow (\Omega^2_{\Bbb C^2} \otimes V)^G$ is surjective.
\end{Theorem}

Then we have the following nice generalized McKay correspondence for
quotient surface singularities: 

\begin{Theorem}\label{Th:Wun}{\rm (\cite{Wunram})}
There is a bijection between the set of special non-trivial
indecomposable reflexive modules $\cal F_i$ and the set of irreducible
components $E_{i}$ via $c_1(\TF_i)E_j = \delta_{ij}$ where $c_1$ is
the first Chern class, and also a one-to-one correspondence with the
set of special non-trivial irreducible representations. 
\end{Theorem} 

As a corollary of this theorem, we get back the original McKay
correspondence for finite subgroups of $SL(2,\Bbb C )$ because in
this case all irreducible representations are special.

\section{$G$-Hilbert schemes and combinatorics}

In this section, we will discuss $G$-Hilbert schemes and a new way to
find the special representations for cyclic quotient singularities by
combinatorics. 

The Hilbert scheme of $n$ points on $\Bbb C^2$ can be described as a
set of ideals: 
$$ \text{Hilb}^n(\Bbb C^2) = \{ I \subset \Bbb C [x,y] \ | \ I\text
{ ideal}, \dim \Bbb C [x,y]/I=n \}. $$
It is a $2n$-dimensional smooth quasi--projective variety. The 
$G$-Hilbert scheme $\text{Hilb}^G(\Bbb C^2)$ was introduced in the
paper 
by Nakamura and the author (\cite{IN})  as follows:
$$ \text{Hilb}^G(\Bbb C^2) = \{ I \subset \Bbb C [x,y] \ | \ I 
~G \text{-invariant ideal, } \Bbb C [x,y]/I\cong \Bbb C [G] \}, $$
where $|G|=n$. This is a union of components of fixed points of
$G$-action on $\text{Hilb}^n(\Bbb C^2)$ and in fact it is just the
minimal resolution of the quotient singularity $\Bbb C^2/G$. It was
proved for $G\subset SL(2,\Bbb C )$ in \cite{IN} first by the properties of
$\text{Hilb}^n(\Bbb C^2)$ and finite group action of $G$ and a McKay 
correspondence in terms of ideals of $G$-Hilbert schemes was stated.

Later Kidoh (\cite{Kidoh}) proved that the $G$-Hilbert scheme for any
small cyclic subgroup of $GL(2,\Bbb C )$ is also the minimal resolution
of the corresponding  cyclic quotient singularities and
Riemenschneider conjectured that the irreducible representations which are given from the ideals of G-Hilbert scheme, so-called Ito-Nakamura type McKay correspondence, are just same as the special representaions which were defined by himself (\cite{Oswald}). Recently A.~Ishii (\cite{Akira}) proved more generally that the
 $G$-Hilbert scheme for any small $G\subset GL(2,\Bbb C )$ is always
  isomorphic to the minimal resolution of the singularity $\Bbb C^2/G$
  and the conjecture is true: 

\begin{Theorem}\label{Th:Akira}{\rm (\cite{Akira})}
Let $G$ be a finite small subgroup of $GL(2,\Bbb C )$.

(i) $G$-Hilbert scheme $\text{Hilb}^G(\Bbb C^2)$ is the minimal
resolution of $\Bbb C^2/G$. 

(ii) For $y\in \text{Hilb}^G(\Bbb C^2)$, denote by $I_y$ the ideal
corresponding to $y$ and let $m$ be the maximal ideal of 
$\cal O_{\Bbb C^2}$ corresponding to the origin $0$. If $y$ is in the
exceptional locus, then, as representations of $G$, we have  
  $$
   I_y/mI_y\cong 
      \begin{cases}
          \rho_i \oplus \rho_0 & \text{ if } y\in E_i \text{ and } 
       y\not\in E_j \text{ for }j \not= i,\\
           \rho_i \oplus \rho_j \oplus \rho_0 &  \text{ if } y\in E_i \cap  E_j,
       \end{cases}
    $$
  where $\rho_i$ is the special representation associated with the
 ireducible exceptional curve $E_i$. 
\end{Theorem}

\begin{Remark}
In dimension two, we can say that the $G$-Hilbert scheme is the same
as a 2-dimensional irreducible component of the  $G$-fixed set of
$\text{Hilb}^n(\Bbb C^2)$. 
A similar statement holds in dimension three, for 
$G\subset SL(3,\Bbb C )$, that is, the $G$-Hilbert scheme is a
3-dimensional irreducible component of the $G$-fixed set of
$\text{Hilb}^n(\Bbb C^3)$, and a crepant resolution of the quotient
singularity $\Bbb C^3/G$ (cf. \cite{Iku}, \cite{BKR}). In this case, note that 
$\text{Hilb}^n(\Bbb C^3)$ is not smooth. 

Moreover, Haiman proved that the $\cal S_n$-Hilbert scheme
$\text{Hilb}^{\cal S_n} (\Bbb C^{2n})$ is a crepant resolution of
$\Bbb C^{2n}/\cal S_n=$ $n$-th symmetric product of $\Bbb C^2$, i.e.,
$$\text{Hilb}^{\cal S_n} (\Bbb C^{2n}) \cong  \text{Hilb}^n(\Bbb C^2)
$$ 
in the process of the proof of $n$! conjecture. (cf. \cite{Hai})
\end{Remark}

From now on, we restrict our considerations to $G\subset GL(2,\Bbb C)$
cyclic. Wunram constructed the generalized McKay correspondence for
cyclic surface singularities in the paper \cite{Wunwun} and we have to
consider the corresponding geometrical informations (the minimal
resolution, reflexive sheaves and so on) to obtain the special
representations.  
Here we would like to give a new characterization of the special
representations in terms of combinatorics. It is much easier to find
the special representation because we don't need any geometrical
objects, but based on the result of $G$-Hilbert schemes.  

Let us discuss the new characterization of the special representations
in terms of combinatorics. 
Let $G$ be the cyclic group $C_{r,a}$, generated by the matrix
$\begin{pmatrix} \epsilon & 0 \\ 0 & \epsilon^a \end{pmatrix}$ where
$\epsilon^r=1$ and gcd$(r,a)=1$ and consider the character map 
$\Bbb C [x,y] \longrightarrow \Bbb C [t]/ {t^r} $ 
given by $x \mapsto t$ and $y \mapsto t^a$. Then we have a
corresponding character for each monomial in $\Bbb C [x,y]$.  

Let $I_p$ be the ideal of the $G$-fixed point $p$ in the $G$-Hilbert
scheme, then we can define the following sets. 

Consider a $G$-invariant subscheme $Z_p\subset \Bbb C^2$ for which
$H^0(Z_p, \cal O_{Z_p})=\cal O_{\Bbb C^2}/I_p$ is the regular
representation of $G$. Then the $G$-Hilbert scheme can be regarded as
a moduli space of such $Z_p$. 

\begin{Definition}
A set $Y(Z_p)$ of monomials in $\Bbb C[x,y]$ is called 
{\it $G$-cluster} if all monomials in $Y(Z_p)$ are not in $I_p$,
and $Y(Z_p)$ can be drawn as a Young diagram of $|G|$ boxes.    
\end{Definition}

\begin{Definition}
For any small cyclic group $G$, let $B(G)$ be the set of monomials
which are not divisible by any $G$-invariant monomial. We call $B(G)$
$G$-basis.
\end{Definition}

\begin{Definition}
If $|G|=r$, then let $L(G)$ be  
$\{ 1, x, \cdots, x^{r-1}, y, \cdots, y^{r-1} \}$, i.e., the set of
monomials which cannot be divided by $x^r$, $y^r$ or $xy$. 
We call it $L$-space for $G$ because the shape of this diagram looks
like the capital letter \lq\lq $L$."
\end{Definition}

\begin{Definition}
The monomial $x^my^n$ is of {\it weight} $k$ if 
$m+an=k$.
\end{Definition}

Let us describe the method to find the special representations of $G$
with these diagrams: 
\begin{Theorem}\label{Th:A}
For a small finite cyclic subgroup of $GL(2, \Bbb C )$, the irreducible
representation $\rho_i$ is special if and only if the corresponding
monomials in $B(G)$ are not contained in the set of monomials
$B(G)\setminus L(G)$. 
\end{Theorem}

\begin{proof}
In Theorem~\ref{Th:RW} (3), we have the definition of the special
representation, and it is not easy to compute all special
representations. However look at the behavior of the monomials in
$\Bbb C[x,y]$ under the map $\Phi_i: {(\Omega^2_{\Bbb C^2})}^G
\otimes (\cal O_{\Bbb C^2} \otimes V_i)^G \rightarrow 
(\Omega^2_{\Bbb C^2} \otimes V_i)^G$ for each representation $\rho_i$: 

First, let us consider the monomial bases of each set.
Let $V_i=\Bbb C e_i$ and $\rho(g)e_i=\e^{-i}$. An element
$f(x,y)dx\wedge dy\otimes \rho_i$ is in 
$(\Omega^2_{\Bbb C^2} \otimes V_i)^G$ 
if and only if
$$ g^*f(x,y)dx\wedge dy\cdot \e^{1+a} \otimes \e^{-i}=f(x,y)dx\wedge dy,$$
that is, 
$$ g^*(f(x,y)dx\wedge dy)= \e^{i-(a+1)}(f(x,y)dx\wedge dy).$$
Therefore the monomial basis for $(\Omega^2_{\Bbb C^2} \otimes V_i)^G$
is the set of monomials $f(x,y)$ such that
$$ g: f(x,y) \mapsto \e^{i-(a+1)} f(x,y)$$
under the action of $G$, that is, monomials of weight $i-(a+1)$.

Similarly, we have the monomial basis for ${(\Omega^2_{\Bbb C^2})}^G$ as
the set of monomials $f(x.y)$ of weight $r-(a+1)$. 

The monomial basis for $(\cal O_{\Bbb C^2} \otimes V_i)^G $ is given 
as the set of monomials $f(x,y)$ of weight  $i$.

Let us check the surjectivity of the map $\Phi_i$. If $\Phi_i$ is
surjective, then the monomial basis in  $(\Omega^2_{\Bbb C^2}
\otimes V_i)^G$ can be obtained as the product of the monomial bases of
two other sets. Therefore the degree of the monomials in
$(\Omega^2_{\Bbb C^2} \otimes V_i)^G$ must be higher than the degree
of the monomials in $(\cal O_{\Bbb C^2} \otimes V_i)^G $.

Now look at the map $\Phi_{a+1}$. The vector space $(\cal O_{\Bbb C^2}
\otimes V_{a+1})^G $ is generated by the monomials of weight $a+1$,
i.e., $x^{a+1}, xy, \cdots, y^{b}$ where $ab =a+1\mod r$. 
On the other hand, $(\Omega^2_{\Bbb C^2} \otimes V_{a+1})^G$ is
generated by the degree $0$ monomial $1$. Then the map $\Phi_{a+1}$ is
not surjective. 

By this, if a monomial of type $x^m y^n$, where $mn\not=0$, is a generator
of $(\cal O_{\Bbb C^2} \otimes V_i)^G $,  
then there exists a monomial $x^{m-1}y^{n-1}$ in $(\Omega^2_{\Bbb C^2}
\otimes V_i)^G$ and the degree become smaller under the map $\Phi_i$.  
This means $\Phi_i$ is not surjective. 

Moreover, if the bases of $(\cal O_{\Bbb C^2} \otimes V_i)^G $ is
generated only by $x^i$ and $y^j$ where $aj\equiv i\mod r$, then  the 
degrees of the monomials in $(\Omega^2_{\Bbb C^2} \otimes V_i)^G$ are
larger and $\Phi_i$ is surjective. Thus we have the assertion.
\end{proof}

\begin{Remark}
From this theorem, we can also say that a representation $\rho_i$ is
special if and only if the number of the generators of the space
$(\cal O_{\Bbb C^2} \otimes V_i)^G $ is 2. 
However, as a module over the 
invariant ring ${\cal O}_{\Bbb C^2}^G$ it is minimally generated by 2 elements. 
In this form, the remark is not new. It follows easily in one direction 
from the remark 
after Theorem 2.1 in Wunram's paper \cite{Wunram}, and in the other direction from Theorem 2.1 in combination with the fact proven in the first appendix of 
that paper that in the case of cyclic quotient surface singularities a 
reflexive module is determined by the "Chern numbers" of its torsionfree 
preimage on the minimal resolution. 
\end{Remark}

\begin{Theorem}\label{Th:B}
Let $p$ be a fixed point by the $G$-action, then we can define an
ideal $I_p$ by the $G$-cluster and the configuration of the
exceptional divisors can be described by these data.
\end{Theorem}

\begin{proof}
The defining equation of the ideal $I_p$ is given by
$$\begin{cases} & x^a=\alpha y^c,\\
          & y^b=\beta x^d, \\
         & x^{a-d} y^{b-c}=\alpha \beta,\end{cases}$$
where $\alpha$ and $\beta$ are complex numbers and 
both $x^a$ and $y^c$ (resp. $y^b$ and $x^d$) correspond to the same 
representation (or character). 

The pair $(\alpha, \beta)$ is a local affine coordinate near the fixed
point $p$ and it is also obtained from the calculation with toric
geometry. Moreover each axis of the affine chart is just a exceptional
curve or the original axis of $\Bbb C^2$. The exceptional curve is
isomorphic to a $\Bbb P^1$ and the points on it is written by the
ratio 
like $(x^a : y^b)$ (resp. $(x^d : y^c)$) which is corresponding to a special representation $\rho_a$ (resp. $\rho_d$). The fixed point $p$ is the intersection point of 2 exceptional curves $E_a$ and $E_d$.

Thus we can get the whole space of exceptional locus by deforming
the point $p$ and patching the affine pieces.   
\end{proof}

We will see a concrete example in the following section. Here we would
like to make one remark as a corollary: 

\begin{Corollary}
For $A_n$-type simple singularities, all $n+1$ affine charts can be
described by $n+1$ Young diagrams of type $(1,\cdots,1,k)$. 
\end{Corollary}

\begin{proof}
In $A_n$ case, $xy$ is always $G$-invariant, hence
$B(G)=L(G)$. Therefore we have $n+1$ $G$-clusters  and each of them
corresponds to the monomial ideal $(x^{k}, y^{n-k+2},xy)$.  
\end{proof} 

\section{Example and related topics}

\def\e{\epsilon}

First, we recall the toric resolution of cyclic quotient singularities
because the quotient space $\Bbb C^2/G$ is a toric variety. 

Let $\Bbb R^2$ be the 2-dimensional real vector space,  
$\{ e^i | i=1,2 \}$ its standard base, $L$ the lattice generated by
$e^1$ and $e^2$,  $N:=L+\sum \Bbb Zv$, where the summation runs over
all the elements $v=1/r(1,a)\in G=C_{r,a}$, and  
$$ \sigma: = \left\{ \sum_{i=1}^2 x_ie^i \in \Bbb R^2, 
\quad x_i \geq 0, \forall i, 1\leq i \leq 2 \right\} $$
the naturally defined  rational convex polyhedral cone in $N_{\Bbb
  R}=N\otimes_\Bbb Z \Bbb R$. The corresponding affine torus embedding
$Y_{\sigma}$ is defined as Spec($\Bbb C[\check{\sigma} \cap M]$),
where $M$ is the dual lattice of $N$ and $\check{\sigma}$ the dual
cone of $\sigma$ in $M_{\Bbb R}$ defined as 
$\check{\sigma} :=\{ \xi\in M_{\Bbb R} | \xi(x) \geq 0 , \forall 
x \in \sigma\} $.  

Then $X=\Bbb C^2/G$ corresponds to the toric variety which is induced
by the cone $\sigma$ within the lattice $N$. 

{\bf Fact 1} We can construct a simplicial decomposition $S$ with  the
vertices on the Newton Boundary, that is, the convex hull of the
lattice points in $\sigma$ except the origin.  
\medskip

{\bf Fact 2}
If $\widetilde X:=X_S$ is the corresponding torus embedding, then
$X_S$ is non-singular. Thus, we obtain the minimal resolution
$\pi=\pi_S:\widetilde X=X_S\longrightarrow \Bbb C^2/G=Y$. 
Moreover, each lattice point of the Newton boundary corresponds to an
exceptional divisor. 

\bigskip

\begin{figure}[htbp]
\begin{center}
\includegraphics{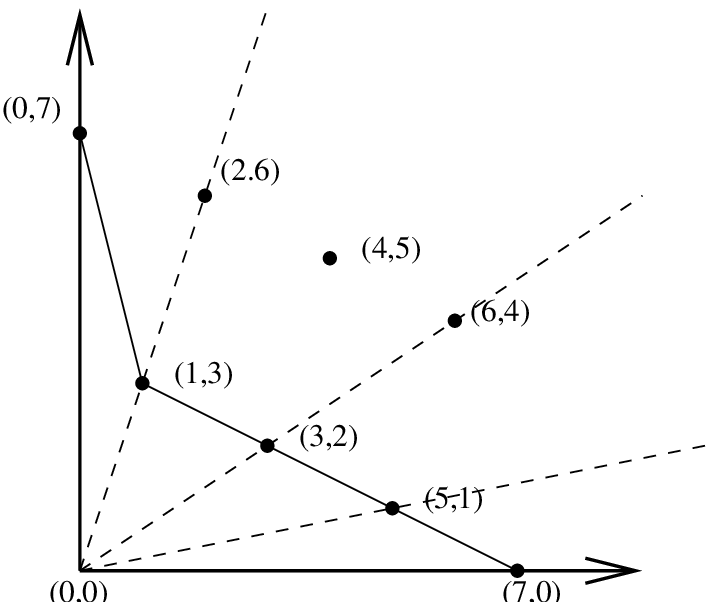}
\caption{toric resolution of $\Bbb C^2/G$}
\label{fig:toric}
\end{center}
\end{figure}

{\bf Example}
Let us look at the example of the cyclic quotient singularity of type
$C_{7,3}$ which is generated by the  matrix  
$\begin{pmatrix} \epsilon & 0 \\ 0 & \epsilon^3 
\end{pmatrix}$ 
where $\epsilon^7=1$. The toric resolution of this quotient
singularity is given by the triangulation of the lattice 
$N\colon=\Bbb Z^2 + \frac{1}{7}(1,3)\Bbb Z$ with the lattice points:
See Figure~\ref{fig:toric}. 

\begin{figure}[htbp]
\begin{center}
\includegraphics{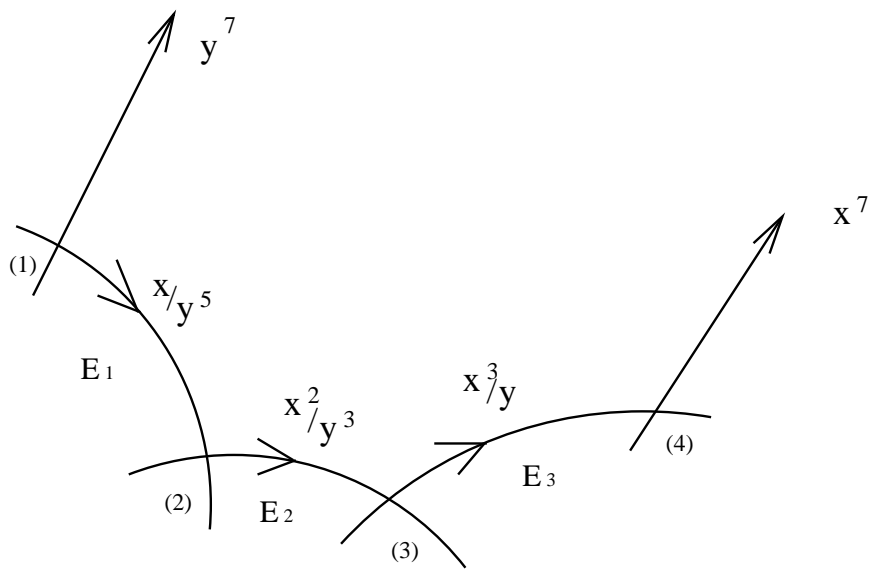}
\caption{configuration of $\tilde X$}
\label{fig:dual}
\end{center}
\end{figure}

From this Newton polytope, we can see that there are 3 exceptional
divisors and the dual graph gives the configuration of the exceptional
components with a deformed coordinate from the original coordinate
$(x,y)$ on $\Bbb C^2$ as in   Figure~\ref{fig:dual}. 

Therefore we have 4 affine pieces in this example and we have 4
coordinate systems corresponding to each affine piece. 
In this picture, we will see the corresponding special irreducible
representations, but we would like to use our method in the previous
section to find the representations. 

\begin{figure}[htbp]
\begin{center}
\includegraphics{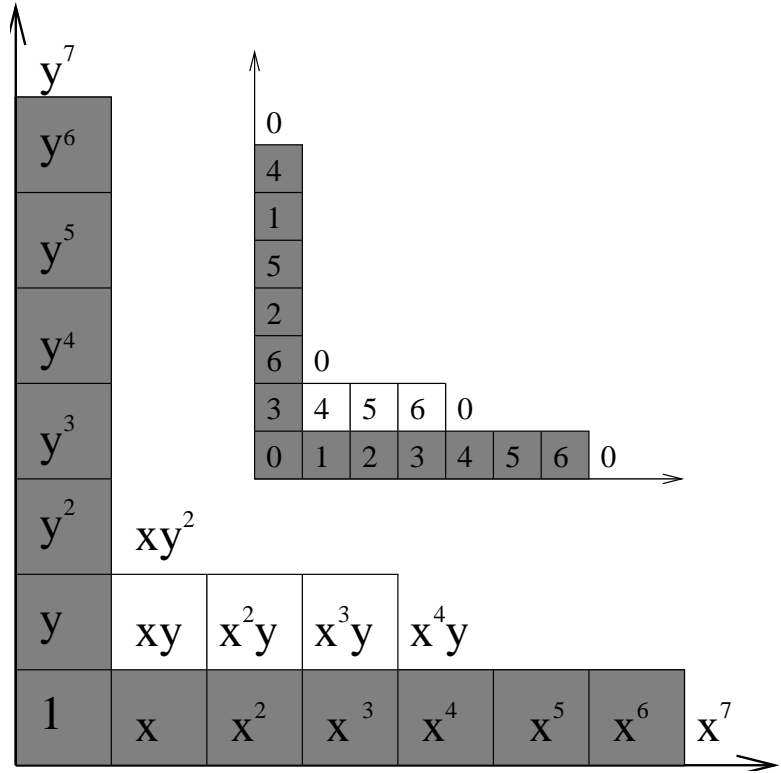}
\caption{$G$-basis $B(G)$ and the characters}
\label{fig:BG}
\end{center}
\end{figure}

Let us draw the diagram which corresponds to the $G$-basis and
$L$-space. First we have the following $G$-basis $B(G)$ and the
corresponding characters in a same diagram. In Figure~\ref{fig:BG} 
we draw the $L$-space as shaded part in $B(G)$. 
  
\begin{figure}[htbp]
\begin{center}
\includegraphics{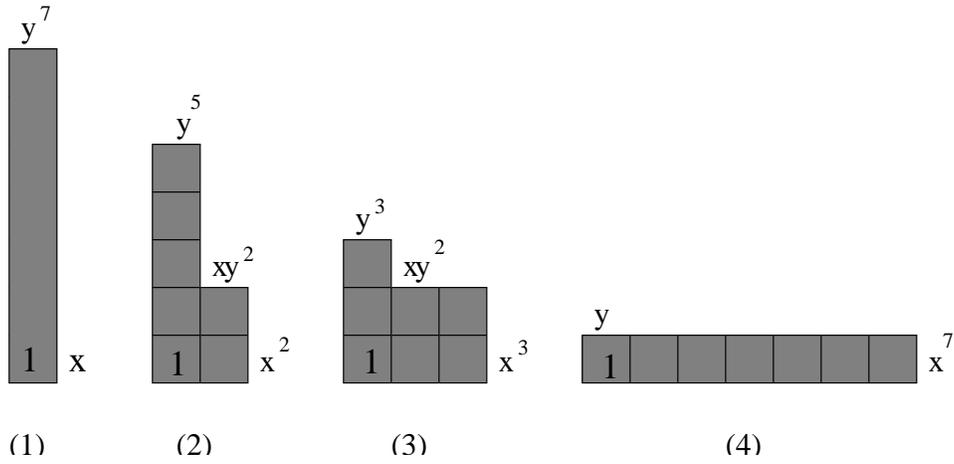}
\caption{$G$-cluster $Y(Z_p)$}
\label{fig:YG}
\end{center}
\end{figure}

Now we have three monomials $xy$, $x^2y$ and $x^3y$ in 
$B(G)\setminus L(G)$ and they correspond to the characters
(resp. representations) 
$4$, $5$ and $6$ (resp. $\rho_4$, $\rho_5$ and $\rho_6$). 
Therefore we can find a set of special representations, that is,
$\{\rho_1,\rho_2,\rho_3\}$, and find the corresponding  
$G$-clusters, representing the origin of the affine charts of the
resolution, can be drawn as 4 Young diagrams and get the corresponding
special representations in this case. See  Figure~\ref{fig:YG}. 

Let us see the meaning of the corresponding $G$-clusters in this
case. From $Y(Z_p)$ for (2), we obtain an ideal 
$I_2=( y^5, x^2, xy^2)$ for the origin of the affine chart (2) in
Figure~\ref{fig:dual}, and the corresponding representations are  
$\rho_1$, $\rho_2$ and $\rho_0$. If we take the  maximal ideal $m$ of
${\cal O}_{\Bbb C^2}$ corresponding to the origin $0$,   
then we have $$I_2/mI_2 \cong \rho_1 \oplus \rho_2 \oplus \rho_0.$$
Similarly we have the ideal $I_3=(y^3, x^3, xy^2)$ and  
$$I_3/mI_3 \cong \rho_2 \oplus \rho_3 \oplus \rho_0.$$
These descriptions coincide with the results of Theorem~\ref{Th:Akira}
for a point at the intersection $E_1 \cap E_2$. 

For any other points $p$  on the exceptional component $E_i$, we must have 
$$I_p/mI_p \cong \rho_i \oplus \rho_0. \qquad \qquad (*)$$
In fact, we can see that a point on the exceptional divisor $E_2$ in this
example was determined by the ratio $x^2:y^3$, that is,  
the corresponding ideal of a point on $E_2$ can be described as
$I_p=(\alpha x^2 - \beta y^3, xy^2-\gamma)$. 
Therefore the ratio $(\alpha : \beta)$ gives the coordinate of the
exceptional curve ($\cong \Bbb P^1$) and we also have $(*)$. 

\bigskip

We discussed special McKay correspondence in 2-dimensional case in
 this paper. In dimension three,  it is convenient to consider crepant
 resolutions as minimal resolutions and we have a much more complicated
 situation. Even in the case $G\subset SL(3,\Bbb C)$, we have
 $H^4(\tilde X, \Bbb Q)\not= 0$ in general. Of course we can use the
 same  definition for the special representations in the higher
 dimensional case, but all non-trivial irreducible representations of
 $G\subset SL(3,\Bbb C)$ are special. On the other hand, the number of
 the exceptional divisors is less than that of the non-trivial
 irreducible representations. 
Therefore, it looks very difficult to generalize this special McKay
 correspondence. That is, we should make a difference, say a kind of
 the grading, in the set of the special (or non-trivial)
 representations like \lq\lq age" of the conjugacy classes.

However, there are good news: In 2000, Craw \cite{Ali} constructed a
cohomological McKay correspondence for the $G$-Hilbert schemes where
$G$ is an abelian group, and in this correspondence we can see the
2-dimensional special McKay correspondence. And recently, the author
found a way to obtain a polytope which corresponds to the
3-dimensional $G$-Hilbert schemes for abelian subgroups in $SL(3,\Bbb
C)$ by combinatorics. There are many crepant resolutions in general in
higher dimension, but the $G$-Hilbert scheme for $G\subset SL(3,\Bbb C)$
is a unique crepant resolution, and the configuration of the
exceptional locus of the special crepant resolution, $G$-Hilbert
scheme, can be determined in terms of a Gr\"obner basis. 
(Let us call this the Gr\"obner method.) Moreover, we can get another
characterization of special representations for cyclic quotient
surface singularities by this Gr\"obner method. So the author is
dreaming of having a more simple and beautiful formulation of the
McKay correspondence in the future.

\end{document}